\documentclass[10pt]{article}  
\usepackage{latexsym,amssymb,amsmath,eepic}

\def\N{\mathbb{N}}

\def\R{\mathbb{R}}

\def\proof{\par\medskip\noindent{\em Proof. }}
\def\eproof{\hfill{$\Box$}\bigskip}

\def\be{\beta}
\def\ga{\gamma}
\def\de{\delta}
\def\si{\sigma}
\def\ep{\varepsilon}

\def\sus{\subset}
\def\ds{\dots}

\newtheorem{thm}{Theorem}
\newtheorem{pro}[thm]{Proposition}

\begin{document}
\title{Polymath's combinatorial proof of the density Hales--Jewett theorem} 
\author{Martin Klazar\footnote{Department of Applied Mathematics, Charles University, Faculty of Mathematics and
Physics, Malostransk\'e n\'am\v est\'\i\ 25, 118 00 Praha, Czech Republic, {\tt klazar@kam.mff.cuni.cz}}}
\maketitle

\begin{abstract}
This is an exposition of the combinatorial proof of the density Hales--Jewett theorem, due to 
D.\,H.\,J. Polymath in 2012. The theorem says that for given $\de>0$ and $k$, 
for every $n>n_0$ every set $A\sus\{1,2,\ds,k\}^n$ with $|A|\ge\de k^n$ contains a combinatorial 
line. It implies Szemer\'edi's theorem, which claims that for given $\de>0$ and $k$, for every $n>n_0$ every set 
$A\sus\{1,2,\ds,n\}$ with $|A|\ge\de n$ contains a $k$-term arithmetic progression.
\end{abstract}

\section{Introduction}

The purpose of this text is to familiarize the author, and possibly the interested reader, with the recent 
remarkable elementary proof \cite{poly,poly_arch} of Polymath (a group of mathematicians, 
see Nielsen \cite{niel} and Gowers \cite{gowe} for more information) for the density Hales--Jewett theorem, 
one of the deepest results in extremal combinatorics/Ramsey theory, which has as an easy corollary the famous theorem 
of Szemer\'edi, indeed the multidimensional generalization thereof. The author hopes 
to use it in his future book on number theory; other similar on-line available fragments are 
\cite{klaz_I,klaz_II,klaz_III}. We begin with recalling the mentioned theorems and introducing some notation. 
Further notation, concepts and auxiliary results will be introduced along the way.

We denote $\N=\{1,2,\ds\}$, $\N_0=\{0,1,\ds\}$ and, for $n\in\N$, $[n]=\{1,2,\ds,n\}$. 
For finite sets $B\ne\emptyset$ and $A$, we call the ratio 
of cardinalities $\frac{|A\cap B|}{|B|}\in[0,1]$ the {\em density of $A$ in $B$} and write $\mu_B(A)$ for it; when $B$ is 
understood from the context, we write just $\mu(A)$ and speak of {\em density of $A$}. Later we consider more general densities.
Densities and the quantities bounding them are denoted by the Greek letters $\mu$, $\de$, $\ep$, $\ga$, $\nu$, $\eta$, $\theta$, $\beta$
and are real numbers from the interval $[0,1]$. A {\em partition} of a set 
$A$ is an expression of $A$ as a disjoint union of possibly empty sets. Note that if $B=\bigcup_{i\in I}B_i$ is a 
partition and $\mu_B(A)\ge\de$, then $\mu_{B_i}(A)\ge\de$ for some $i$. For $a,d,k\in\N$, the $k$-element set 
$$
\{a,a+d,a+2d,\ds,a+(k-1)d\}
$$
is the {\em $k$-term arithmetic progression.} The following is the famous {\em Szemer\'edi's theorem} \cite{szem}.

\begin{thm}\label{szem}
For every $\de>0$ and $k\in\N$, there exists an $n_0\in\N$ such that for every $n>n_0$ every set
$A\sus[n]$ with $\mu(A)\ge\de$ contains a $k$-term arithmetic progression.
\end{thm}
Precursor of Szemer\'edi's theorem
was its color version, the {\em van der Waerden theorem} \cite{vdWa} asserting that for every $r,k\in\N$, for
any $n>n_0$ in any partition $[n]=A_1\cup A_2\cup\ds\cup A_r$ a block $A_i$ contains a $k$-term 
arithmetic progression. Clearly, Szemer\'edi's theorem implies van der Waerden theorem.

For $k,n\in\N$, the set $[k]^n$ consists of all $k^n$ $n$-tuples, called {\em words}, $x=(x_1,x_2,\ds,x_n)$ with $x_i\in[k]$.
For $x\in[k+1]^n\backslash[k]^n$ and $i\in[k]$, we denote by $x(i)$ the word obtained from $x$ by replacing each 
occurrence of $k+1$ by $i$. The $k$-element subset of $[k]^n$ of these words,
$$
L(x)=\{x(i)\;|\;i\in[k]\},
$$
is the {\em combinatorial line (determined by $x$)}. In 1963 Hales and Jewett \cite{hale_jewe} proved that for 
every $r,k\in\N$, for any $n>n_0$ in any partition $[k]^n=A_1\cup A_2\cup\ds\cup A_r$ a block $A_i$ contains 
a combinatorial line. The stronger density version of this theorem was achieved by Furstenberg and 
Katznelson in 1991 \cite{furs_katz} (they proved the special case $k=3$ earlier in 
\cite{furs_katz_89}) by ergodic methods, developed by Furstenberg \cite{furs} in his proof of Szemer\'edi's theorem.
Thus, the {\em density Hales--Jewett theorem} asserts the following.

\begin{thm}\label{denshalejewe}
For every $\de>0$ and $k\in\N$, there exists an $n_0\in\N$ such that for every $n>n_0$ every set
$A\sus[k]^n$ with $\mu(A)\ge\de$ contains a combinatorial line.
\end{thm}
We shall prove Theorem~\ref{denshalejewe}, following Polymath's proof in \cite{poly}. Theorem~\ref{denshalejewe} 
implies Theorem~\ref{szem}, with the same $k$, by means of the bijection
$$
f:\;[k]^n\to[k^n],\ f(x)=f((x_1,x_2,\ds,x_n))=1+\sum_{i=1}^n(x_i-1)k^{i-1}
$$
which sends combinatorial lines to $k$-term arithmetic progressions and, being bijection, preserves densities; for 
the color versions of the theorems the simpler mapping $x\mapsto x_1+x_2+\ds+x_n$ suffices for the reduction.

{\em Multidimensional Szemer\'edi's theorem} claims that for every $\de>0$, $r\in\N$ and finite set $H\sus\N^r$, there
exists an $n_0\in\N$ such that for every $n>n_0$ every set $A\sus[n]^r$ with $\mu(A)\ge\de$ contains 
a copy of $H$ of the form
$a+dH$, $a\in\N^r,d\in\N$. The particular case with $r=2$ and $H=\{(1,1),(1,2),(2,1)\}$ is the {\em corner 
theorem} which was derived by Ajtai and Szemer\'edi \cite{ajta_szem} from Szemer\'edi's theorem. As explained in
\cite{poly} and \cite{gowe}, the proof of Theorem~\ref{denshalejewe} in \cite{poly} is inspired by and modelled after 
the increment density argument in \cite{ajta_szem}.  

\section{The proof of Theorem~\ref{denshalejewe}}

The {\em combinatorial subspace $S$ of $[k]^n$ with dimension $d$}, $d\le n$, determined by the word $x\in[k+d]^n$ 
such that each letter $k+1,k+2,\ds,k+d$ appears in $x$ at least once, is the $k^d$-element subset of $[k]^n$
$$
S=S(x)=\{x(y)\;|\;y\in[k]^d\}
$$
where $x(y)$ is the word obtained from $x$ by replacing each occurrence of $k+i$ by $y_i$, $i=1,2,\ds,d$.
In other words, $S\sus[k]^n$ is a $d$-dimensional combinatorial subspace of $[k]^n$ if and only if there exist a word
$z\in[k]^n$ and $d$ nonempty and disjoint subsets $X_l\sus[n]$, $1\le l\le d$, such that 
$$
x\in S\iff\mbox{$x_i=z_i$ if $i\in[n]\backslash(X_1\cup\dots\cup X_d)$ and $x_i=x_j$ if $i,j\in X_l$.}
$$
The elements of $[n]$ in the union $X_1\cup\ds\cup X_d$ are the {\em free coordinates of $S$} and those not in it 
are the {\em fixed coordinates of $S$}.
The $1$-dimensional combinatorial subspaces are exactly combinatorial lines. From now we omit for brevity 
`combinatorial' for subspaces and lines. 
The words $[k+d]^n\backslash\bigcup_{l=k+1}^{k+d}([k+d]\backslash\{l\})^n$ and $d$-dimensional subspaces
of $[k]^n$ correspond via the mapping $x\mapsto S(x)$, and this is a $d!$-to-one correspondence as $S(x)=S(x')$ iff
$x$ and $x'$ can be identified by permuting the letters $k+1,k+2,\ds,k+d$. The set of words $[k]^d$ and any 
$d$-dimensional subspace $S(x)\sus[k]^n$ are in bijection via $y\mapsto x(y)$. This bijection sends 
the $e$-dimensional subspaces of $[k]^d$, $e\le d$, to the $e$-dimensional subspaces of $[k]^n$ contained in $S(x)$,  
and this is in fact a bijection.

We capture the {\em density increment argument} by the next proposition. 
\begin{pro}\label{increment}
There is a function 
$$
c=c(k,\de):\;\N\times(0,1)\to(0,1), 
$$ 
nondecreasing in $\de$ for every $k$, such that for every $k,d\in\N$ and $\de\in(0,1)$, there is an $n_0$ such that for 
every $n>n_0$ and every set $A\sus[k]^n$ with $\mu(A)\ge\de$ and containing no line, there exists a subspace $S\sus[k]^n$ 
with dimension $d$ and 
$$
\mu_S(A)\ge\mu(A)+c\ge\de+c.
$$
\end{pro}
We fix $k\ge 2$ and derive Theorem~\ref{denshalejewe} from Proposition~\ref{increment}. Suppose $\de>0$ is given and 
let $c=c(k,\de)>0$. By Proposition~\ref{increment}, for $d=1$ there is an $n_0$ such that if $n>n_0$, 
$A\sus[k]^n$ has $\mu(A)\ge\de$ and avoids lines, then we get (by the bijection between 
$S$ and $[k]^d$) a set $A'\sus[k]^d=[k]$ that has $\mu(A')\ge\de+c$. For 
$d=n_0+1$ we have the conclusion for every $n>n_1$ for some $n_1$ and get $A'\sus[k]^d=[k]^{n_0+1}$ 
free of lines and with $\mu(A')\ge\de+c$. We apply to $A'\sus[k]^{n_0+1}$ Proposition~\ref{increment} again
and get $A''\sus[k]$ with $\mu(A'')\ge(\de+c)+c=\de+2c$. We iterate the process and define 
inductively in a clear way numbers $n_2,n_3,\ds,n_t$ where $t=\lfloor 1/c\rfloor$. 
For $n>n_t$, every set $A\sus[k]^n$ with $\mu(A)\ge\de$ contains a line, for else repeated applications 
of Proposition~\ref{increment} produce at the end a subset of $[k]$ with density at least $\de+(t+1)c>1$, 
which cannot exist.

The density increment $c$ of Proposition~\ref{increment} arises in two steps, embodied in the next 
two propositions. For $k\ge2$ and $i\in[k-1]$, a set $D\sus[k]^n$ is called {\em $ik$-insensitive} if 
$x\in D\Rightarrow x'\in D$ for any word $x'$ obtained from $x$ by changing some occurrences of $k$ to $i$ 
or vice versa. If $D\sus S\sus[k]^n$, where $S$ is a $d$-dimensional subspace, we say that 
$D$ is {\em $ik$-insensitive in $S$} if $D'\sus[k]^d$ is $ik$-insensitive where $D'$ is the image of $D$ in the 
bijection between $S$ and $[k]^d$. A set $D\sus[k]^n$ is a {\em $k$-set} if $D=\bigcap_{i=1}^{k-1}D_i$ where 
each $D_i\sus[k]^n$ is $ik$-insensitive. We define $k$-sets $D\sus S$ in a $d$-dimensional subspace $S\sus[k]^n$ 
similarly, via the bijection between $S$ and $[k]^d$. In the next two propositions we may assume that $k\ge3$ since 
they will be used for such $k$.

\begin{pro}\label{partition}
Let $k,d\in\N$ and $\ep>0$. There exists an $n_0$ such that for every $n>n_0$, every $k$-set $D\sus[k]^n$ has a
partition 
$$
D=S_1\cup\ds\cup S_t\cup F
$$
into $d$-dimensional subspaces $S_i\sus[k]^n$ and a set $F\sus[k]^n$ with $\mu(F)<\ep$.
\end{pro}

\begin{pro}\label{inssens}
There is a function
$$
\ga=\ga(k,\de):\;\N\times(0,1)\to(0,1),
$$ 
nondecreasing in $\de$ for every $k$, such that for every $k,r\in\N$ and $\de\in(0,1)$, there is an $n_0$ such that
for every $n>n_0$ and every set $A\sus[k]^n$ with $\mu(A)\ge\de$ and containing no line, there exist an 
$r$-dimensional subspace $W\sus[k]^n$ and a $k$-set $D\sus W$ in $W$ satisfying
$$
\mu_W(D)\ge\ga\ \mbox{ and }\ \mu_D(A)\ge\mu(A)+\ga\ge\de+\ga.
$$ 
\end{pro}
We fix $k\ge 2$ and derive Proposition~\ref{increment} from Propositions~\ref{inssens} and \ref{partition}. 
Let $d$ and $\de$ 
be given. We set $\ga=\ga(k,\de)>0$ and take the $n_0$ of Proposition~\ref{partition} corresponding to $k,d$ 
and $\ep=\ga^2/2$. Then we take $n_1$ such that for $n>n_1$ Proposition~\ref{inssens} holds with $k,r=n_0+1$ and $\de$.
Now let $n>n_1$ and suppose a set $A\sus[k]^n$ with $\mu(A)\ge\de$ and free of lines is given. 
There exist a subspace $W$ and a $k$-set $D\sus W$ of Proposition~\ref{inssens} such that $\mu_W(D)\ge\ga$, $\mu_D(A)\ge\mu(A)+\ga$ and $W$ has dimension $n_0+1$. 
Thus $D$ partitions as in Proposition~\ref{partition}, with $[k]^{n_0+1}$ corresponding to $W$ in place of 
the $[k]^n$ in Proposition~\ref{partition} and $\ep=\ga^2/2$. Let $D=E\cup F$ be a partition where $E$ is a disjoint 
union of $d$-dimensional subspaces of $W$ and $\mu_W(F)<\ep$. Since $\mu_D(A)\ge\mu(A)+\ga$, 
$\mu_D(F)=\mu_W(F)/\mu_W(D)<\ep/\ga=\ga/2$ and $\mu_F(A)\le 1$, we get $\mu_E(A)\ge\mu(A)+\ga/2$. By averaging, 
there is a $d$-dimensional subspace $S$ of $W$ (contained in $E$) with $\mu_S(A)\ge\mu(A)+\ga/2$. 
Proposition~\ref{increment} follows, with $c(k,\de)=\ga(k,\de)/2$. 

Thus to prove Theorem~\ref{denshalejewe} it suffices to deduce Propositions~\ref{inssens} and \ref{partition}. We 
shall proceed by induction on $k$. We start by proving Theorem~\ref{denshalejewe} for $k=2$ and then for every 
$k\ge3$ derive Propositions~\ref{inssens} and \ref{partition} from validity of Theorem~\ref{denshalejewe} for $k-1$. 
The derivations rely on formally stronger but equivalent forms of Theorem~\ref{denshalejewe}, Propositions~\ref{multidhj}
and \ref{pravddhjthm}. We get the implications
$$
\mathrm{T2}_2\Rightarrow\mathrm{P4}_3\;\&\;\mathrm{P5}_3\Rightarrow\mathrm{P3}_3
\Rightarrow\mathrm{T2}_3\Rightarrow\mathrm{P4}_4\;\&\;\mathrm{P5}_4\Rightarrow\mathrm{P3}_4
\Rightarrow\mathrm{T2}_4\Rightarrow\ds,
$$
which establish Theorem~\ref{denshalejewe} for every $k\ge 2$. We start with the easy case $k=2$ and then 
prepare some results for the derivation of Propositions~\ref{inssens} and \ref{partition}.

The words of $[2]^n$ 1-1 correspond to the subsets $X\sus[n]$, and the lines 1-1 correspond
to the {\em inclusion pairs:} pairs $X\sus Y\sus[n]$ with $Y\ne X$. Thus Theorem~\ref{denshalejewe} for $k=2$ follows from the 
next classical {\em Sperner's theorem} \cite{sper}.

\begin{pro}\label{sperner}
If $F$ is a family of subsets of $[n]$ containing no inclusion pair (i.e., $F$ is an antichain to $\sus$) then 
$$
|F|\le\binom{n}{\lfloor n/2\rfloor}\ll\frac{2^n}{\sqrt{n}}.
$$
\end{pro}
\proof
Let $F$ be an antichain of subsets of $[n]$. The maximal chains $\{X_0\sus X_1\sus\ds\sus X_n=[n]\}$, $|X_i|=i$, 
1-1 correspond to the $n!$ permutations $\pi$ of $[n]$ via $\pi\mapsto C_{\pi}=\{\emptyset,\{\pi(1)\},\{\pi(1),\pi(2)\},
\ds,\pi([n])=[n]\}$. We double count the pairs $(\pi,X)$ such that $X\in C_{\pi}\cap F$. Grouping the pairs by $\pi$
we get that their number is $\le n!$ as $|F\cap C_{\pi}|\le1$ for each $\pi$. Grouping them by $X$ we get that
their number is exactly $\sum_{X\in F}|X|!(n-|X|)!$, since the summand equals to the number of $\pi$ with $X\in C_{\pi}$. 
Hence
$$
\sum_{X\in F}|X|!(n-|X|)!\le n!.
$$
Since $\binom{n}{j}\le\binom{n}{\lfloor n/2\rfloor}$ for any $0\le j\le n$, 
$\lfloor n/2\rfloor!(n-\lfloor n/2\rfloor)!\le|X|!(n-|X|!)$ and dividing by $\lfloor n/2\rfloor!(n-\lfloor n/2\rfloor)!$
yields the stated inequality.
\eproof

One may generalize Theorem~\ref{denshalejewe} to subspaces but this is not 
really stronger than the original theorem.

\begin{pro}\label{multidhj}
Let $k\in\N$, $k\ge2$, be given. Assuming Theorem~\ref{denshalejewe} for $k$, it follows that for every 
$\de>0$ and $d\in\N$, there exists an $n_0\in\N$ such that for every $n>n_0$ every set $A\sus[k]^n$ with $\mu(A)\ge\de$ 
contains a $d$-dimensional subspace.
\end{pro}
\proof
We proceed by induction on $d$ where the case $d=1$ is Theorem~\ref{denshalejewe}. Now suppose that $d\ge2$ and 
the result holds for $d-1$ (and every $\de$). Observe that if $n=n_1+n_2$, $n_i\in\N$, and $A\sus[k]^n$ with
$\mu(A)\ge\de$, then
$$
\mu(A_1)\ge\de/2\ \mbox{ for }\ A_1=\{x\in[k]^{n_1}\;|\;\mu(\{y\in[k]^{n_2}\;|\;(x,y)\in A\})\ge\de/2\}
$$
(interpreting $(x,y)$ in the obvious way as an element of $[k]^n$). Let $\de>0$ be given.
We take an $n_2$ such that the result holds (with $n=n_2$)
for $d-1$ and density $\de/2$ and then take an $n_1$ such that the conclusion of Theorem~\ref{denshalejewe} holds 
for every $n>n_1$, with density $\de/2(k+d-1)^{n_2}$. Suppose that $n>n_1+n_2$ and $A\sus[k]^n$ has $\mu(A)\ge\de$.
Then, using the observation, inductive assumption and pigeonhole principle, we get a set $A_1\sus[k]^{n-n_2}$
with $\mu(A_1)\ge\de/2(k+d-1)^{n_2}$ and a $(d-1)$-dimensional subspace $S\sus[k]^{n_2}$ such that
$(x,y)\in A$ for every $x\in A_1$ and $y\in S$. By Theorem~\ref{denshalejewe}, $A_1$ contains a line $L$.
Hence $\{(x,y)\;|\;x\in L,y\in S\}$ is the desired $d$-dimensional subspace contained in $A$.
\eproof

We will use the fact that almost all words in $[k]^n$ have almost precisely $n/k$ occurrences of each of the $k$ letters.
\begin{pro}\label{cernov}
Let $k,n\in\N$, $j\in[k]$ and $A\sus[k]^n$ be the set of words with the number of occurrences of $j$ outside 
the interval $[n/k-n^{2/3},n/k+n^{2/3}]$. Then $\mu(A)<n^{-1/3}$.
\end{pro}
\proof
(We say more on the tools used in  Subsection~\ref{auxres}.)
For $i\in[n]$ and $x\in[k]^n$, let $f_i(x)=1$ if $x_i=j$ and $f_i(x)=0$ else. Then the function 
$f=f_1+\ds+f_n$ counts occurrences of $j$ in $x$, has mean $P=n/k$ (sum of the means of the $f_i$) and 
variance $V=(n/k)(1-1/k)<n$ ($V$ is the mean of $f^2$ minus the square of the mean of $f$, which by linearity of means 
and independence of the $f_i$ gives $n/k+n(n-1)/k^2-(n/k)^2$).
By \v{C}eby\v{s}ev's inequality, $\mu(\{x\in[k]^n\;|\;|f(x)-P|>\lambda\sqrt{V}\})<\lambda^{-2}$ 
for any $\lambda>0$. Setting $\lambda=n^{1/6}$ gives the result.
\eproof

For $k,n\in\N$, we have on $[k]^n$ the uniform density $\mu$, given by 
$\mu(\{x\})=1/k^n$. For $k\ge2$ and the parameter $m\le n$, we define another, non-uniform, density $\mu'_m$ on $[k]^n$ by 
$$
\mu'_m(\{x\})=\frac{|\{(J,y,z)\in M\;|\;(J,y,z)=x\}|}{|M|},
$$
for
$$
M=\{(J,y,z)\;|\;J\sus[n],|J|=m,y\in[k-1]^J,z\in[k]^{[n]\backslash J}\}
$$
($A^B$ denotes, for sets $A$ and $B$, the set of all mappings from $B$ to $A$), where any triple 
$(J,y,z)$ in $M$ is interpreted as $x\in[k]^n$ by setting $x_i=y_i$ if $i\in J$ and $x_i=z_i$ else.
(We say more about densities in Subsection~\ref{auxres}.)

\begin{pro}\label{mucarka}
If $\eta\in(0,1)$ and $k,m,n\in\N$ satisfy $k\ge2$, $n\ge(12k/\eta)^{12}$ and $m\le n^{1/4}$, then 
for every set $A\sus[k]^n$ we have
$$
|\mu'_m(A)-\mu(A)|<\eta.
$$
\end{pro}
\proof
This is a particular case of the more general Proposition~\ref{prumdens}, which we 
prove later. 
\eproof

\noindent
To deduce Proposition~\ref{partition}, we need Propositions~\ref{multidhj} and \ref{mucarka}.

\subsection{Derivation of Proposition~\ref{partition}}

In this subsection we fix a $k\in\N$, $k\ge3$, assume that Theorem~\ref{denshalejewe} holds for $k-1$ (and every 
$\de>0$) and deduce from this Proposition~\ref{partition} for $k$. The main step is to get the required partition if 
$D$ is an $ik$-insensitive set for an $i\in[k-1]$; the full proposition follows inductively by iteration. 
We may of course set $i=1$.

\begin{pro}\label{partition1}
For every $d\in\N$ and $\ep>0$, there exists an $n_0$ such that for every $n>n_0$ every $1k$-insensitive set
$D\sus[k]^n$ has a partition
$$
D=S_1\cup\ds\cup S_t\cup F
$$ 
into $d$-dimensional subspaces $S_i\sus[k]^n$ and a set $F\sus[k]^n$ with $\mu(F)<\ep$.
\end{pro}
\proof
Let $d\in\N$ and $\ep>0$ be given. Applying Theorem~\ref{denshalejewe} for $k-1$ and Proposition~\ref{multidhj},  
we take an $m\in\N$, $m\ge d$, 
such that every set $A\sus[k-1]^m$ with $\mu(A)\ge\ep/3$ contains a $d$-dimensional subspace. We set  
$$
n_0=\lceil3\ep^{-1}mk^{m-d}(k+d)^m+m^4+(36k/\ep)^{12}\rceil. 
$$
Let $n>n_0$ and $D\sus[k]^n$ be a $1k$-insensitive set. We may assume that 
$\mu(D)\ge\ep$ for else we set at once $F=D$. We construct, for $r=0,1,\ds$, sets 
$D=D_0\supset D_1\supset\ds\supset D_r$ and $\emptyset=J_0\sus J_1\sus\ds\sus J_r\sus[n]$, $|J_j|=jm$, with 
the properties that (i) for each $x\in[k]^{J_r}$, the set 
$$
(D_r)_x=\{y\in[k]^{[n]\backslash J_r}\;|\;(x,y)\in D_r\}
$$ 
is $1k$-insensitive, (ii) $D\backslash D_r$ partitions into $d$-dimensional subspaces and (iii) 
$\mu(D_j\backslash D_{j+1})\ge\ep k^{d-m}(k+d)^{-m}/3$ for $j=0,1,\ds,r-1$. Such sets trivially 
exist if $r=0$, namely $D_r=D$ and  $J_r=\emptyset$. We claim that as long as $\mu(D_r)\ge\ep$, the construction can be 
continued. This establishes
the proposition: since $r\le3\ep^{-1}k^{m-d}(k+d)^m$ (by (iii)), 
the construction has to terminate for some $r$ (by the definition of $n_0$, $n$ is so large that without terminating
we hit the contradiction $\mu(D)>1$), and then $\mu(D_r)<\ep$ and $D\backslash D_r$ 
partitions into $d$-dimensional subspaces.

To prove the claim we assume that $\mu(D_r)\ge\ep$, which is true if $r=0$. In the initial step when $r=0$ and 
$J_r=\emptyset$, we modify the following construction, which is described for the general step, accordingly by omitting 
the $x$-coordinate. The (uniform) average of the values $\mu((D_r)_x)$, taken over 
all $x\in[k]^{J_r}$, equals $\mu(D_r)$ and so is at least $\ep$. Hence the same average of $\mu_m'((D_r)_x)$, where 
$\mu'_m$ is the (non-uniform) density on $[k]^{[n]\backslash J_r}$ introduced before Proposition~\ref{mucarka}, 
is at least $\ep-\eta=2\ep/3$, due to Proposition~\ref{mucarka} with $\eta=\ep/3$. In other words, 
density of the subset of the quadruples $(x,J,y,z)$, where $x\in[k]^{J_r}$, $J\sus[n]\backslash J_r$, $|J|=m$, $y\in[k-1]^J$
and $z\in[k]^{[n]\backslash(J_r\cup J)}$, satisfying $(x,J,y,z)\in D_r$, in the set of all quadruples, is 
at least $2\ep/3$. Hence there is a $J$ such that density of the triples $(x,y,z)$ (from the stated 
domains) with $(x,J,y,z)\in D_r$ is at least $2\ep/3$. And for this fixed $J$, density of the pairs $(x,z)$ for which 
$\mu(\{y\in[k-1]^J\;|\;(x,J,y,z)\in D_r\})\ge\ep/3$ is at least $\ep/3$. We set 
$$
J_{r+1}=J_r\cup J. 
$$
By the choice of $m$, each of these sets of words $y$ contains a $d$-dimensional subspace $U'_{x,z}\sus[k-1]^J$.
Since $(D_r)_x$ is $1k$-insensitive, for this $J$ and for each of these pairs $(x,z)$ there is a 
$d$-dimensional subspace $U_{x,z}\sus[k]^J$ such that $(x,y,z)\in D_r$ for every $y\in U_{x,z}$. 
By the pigeon-hole principle, there is a single $d$-dimensional subspace $U\sus[k]^J$ such that the set
$$
T=\{(x,z)\in[k]^{J_r}\times[k]^{[n]\backslash J_{r+1}}\;|\;(x,y,z)\in D_r\mbox{ for every $y\in U$}\}
$$
has density at least $\ep/3(k+d)^m$. Note that for each $x$ the set of $z$ with $(x,z)\in T$ is $1k$-insensitive. 
We set
$$
D_{r+1}=D_r\backslash(T\times U)
$$
where $T\times U$ means the words in $[k]^n$ that restrict, for some $(x,z)\in T$ and $y\in U$, on $J_r$ to $x$, 
on $[n]\backslash J_{r+1}$ to $z$ and on $J$ to $y$. Clearly, (ii) holds because $T\times U=D_r\backslash D_{r+1}$ is a disjoint 
union of $d$-dimensional subspaces. Property (i) holds too because for every $x\in[k]^{J_r}$ and 
$y\in[k]^J$, the set of $z\in[k]^{[n]\backslash J_{r+1}}$
with $(x,y,z)\in D_{r+1}$ is $1k$-insensitive. (Consider $u=(x,y,z)\in D_{r+1}$ and $u'=(x,y,z')$
in which $z'$ arises from $z$ by some exchanges of $1$s and $k$s. Then $u'\in D_r$ by the $1k$-insensitivity of 
$(D_r)_x$. If $u'\in T\times U$ then $y\in U$ and $(x,z')\in T$, hence $(x,z)\in T$ as noted above, 
and $u\in T\times U$, which is not the case. So $u'\not\in T\times U$ and  $u'\in D_{r+1}$.) Finally, the density of $T\times U$ in $[k]^n$
equals to the density of $T$ in $[k]^{J_r}\times[k]^{[n]\backslash J_{r+1}}$ times the density of $U$ in $[k]^J$, which 
is at least $\ep/(3(k+d)^mk^{m-d})$. Thus $J_{r+1}$ and $D_{r+1}$ have the required properties (i)--(iii).
\eproof

\noindent
{\em Remark.} The decrease of density of $T\times U$ compared to $T$, caused by density of $U$, 
seems to be overlooked by Polymath---they claim \cite[bottom of p. 1320]{poly} that $T\times U$ has density at 
least $\eta(k+d)^{-m}$ (i.e., $\ep/3(k+d)^m$ in our notation), which reflects in the statement of 
\cite[Lemma 8.1]{poly}. 

We prove Proposition~\ref{partition}. We proceed by induction on the size of intersection defining $D$. 
Let $j\in[k-1]$. We assume that for every $d$ and $\ep>0$ Proposition~\ref{partition} holds for all sets
of the form $D=D_1\cap D_2\cap\ds\cap D_j$ where $D_i\sus[k]^n$ is $ik$-insensitive; for $j=1$ this is true
by Proposition~\ref{partition1}. From this we deduce (if $j<k-1$) that Proposition~\ref{partition} holds 
for all sets $D$ corresponding to the increased parameter value $j+1$. For $j=k-1$ we get 
the original Proposition~\ref{partition}. 

So let $d$ and $\ep>0$ be given. We take $n_0$ such that for every $n>n_0$ our inductive assumption (for $j<k-1$) 
holds for subspaces dimension $d$ and bound on the density of the residual set $\ep/2$. Then we take $n_1$ such that 
for every $n>n_1$ the conclusion of Proposition~\ref{partition1} holds for subspaces dimension $n_0+1$ and 
bound on the density of the residual set $\ep/2$. Now suppose that $n>n_1$ and $D=D_1\cap\ds\cap D_{j+1}$
where $D_i\sus[k]^n$ is $ik$-insensitive. Using Proposition~\ref{partition1} we obtain a partition 
$D_{j+1}=T_1\cup\ds\cup T_s\cup F$ such that the $T_i$ are $(n_0+1)$-dimensional subspaces and
$\mu(F)<\ep/2$. We have
$$
D=\bigcap_{i=1}^{j+1}D_i=\bigcup_{i=1}^s(T_i\cap D_1\cap\ds\cap D_j)\cup(F\cap D_1\cap\ds\cap D_j).
$$
Clearly, each $T_i\cap D_h$ is $hk$-insensitive in $T_i$ and thus using the inductive assumption 
for $j$ we can express (working in $[k]^{n_0+1}$ via the bijection with $T_i$) each set 
$T_i\cap D_1\cap\ds\cap D_j=(T_i\cap D_1)\cap\ds\cap(T_i\cap D_j)$ 
as a disjoint union of $d$-dimensional subspaces (in $T_i$ and thus in $[k]^n$) and a residual 
set $F_i$ with $\mu_{T_i}(F_i)<\ep/2$. These subspaces, taken for all $i=1,2,\ds,s$, and the set
$E=F_1\cup\ds\cup F_s\cup(F\cap D_1\cap\ds\cap D_j)$ form the desired partition of $D$ because
$\mu(E)\le\mu(F_1\cup\ds\cup F_s)+\mu(F)<\max_i\mu_{T_i}(F_i)+\ep/2<\ep$. This concludes the derivation of 
Proposition~\ref{partition}.

\subsection{The equal-slices densities $\nu$ and $\tilde{\nu}$}\label{auxres}

We move to the second and more complicated half of the proof of Theorem~\ref{denshalejewe}, 
the derivation of Proposition~\ref{inssens} for $k$ from Theorem~\ref{denshalejewe} for $k-1$. 
Similarly to the role of  Proposition~\ref{multidhj} in the first half of the proof, 
we need a stronger version of Theorem~\ref{denshalejewe}, Proposition~\ref{pravddhjthm}, which says that any positively dense 
set $A\sus[k]^n$ contains, for large $n$, a set of lines with positive density. However, Proposition~\ref{cernov} shows
that this cannot hold for the uniform density. Consider the set $A\sus[2]^n$ of words in which the numbers of occurrences 
of $1$ and $2$ deviate from $n/2$ by less than $n^{2/3}$. Then $\mu(A)\to 1$ as $n\to\infty$ but at the same time 
$\mu(M)\to0$ for the set $M\sus[3]^n$ of lines contained in $A$, because the inclusion $L(x)\sus A$, 
$x\in[3]^n$, forces $x$ to have at most $2n^{2/3}$ occurrences of $3$, and such $x$ have in $[3]^n$ density going to 
$0$. Fortunately, the strengthening  holds for a non-uniform density, the equal-slices density $\nu$ that we define 
in a moment, and one can go from the uniform to equal-slices density and back. Since $\nu$ does not behave well to
restrictions to subspaces, we need to work also with a variant density $\tilde{\nu}$ fixing this problem, which for large $n$ 
differs from $\nu$ only little. We begin with discussing densities in general and then introduce the densities $\nu$ and 
$\tilde{\nu}$.

A {\em density} on a finite set $B\ne\emptyset$ is a mapping $\mu'$ from the set of all subsets of $B$ to 
the interval $[0,1]$, such that
$\mu'(B)=1$ and $\mu'(A\cup A')=\mu'(A)+\mu'(A')$ whenever $A,A'\sus B$ and $A\cap A'=\emptyset$. Thus $\mu'(\emptyset)=0$
and $\mu'$ is uniquely determined by its values on singletons. Any choice of values $\mu'(\{x\})\ge0$, $x\in B$, with 
$\sum_{x\in B}\mu'(\{x\})=1$ gives a density: $\mu'(A)=\sum_{x\in A}\mu'(\{x\})$ for any $A\sus B$. We have been using
the {\em uniform density $\mu$}, defined by $\mu(\{x\})=1/|B|$ for any $x\in B$, and before Proposition~\ref{mucarka} we 
met the non-uniform density $\mu'_m$. We reserve the letter $\mu$ for the uniform density and primed 
$\mu'$ for general, possibly non-uniform, density. 

Suppose $B$ is a finite set with a density $\mu'$. If $f:\;B\to\R$, the {\em average, or mean, of the function $f$ 
(with respect to $\mu'$)} is
$$
\sum_{x\in B}f(x)\cdot\mu'(\{x\}). 
$$
We recall a few useful properties of averages, which we already used in the proof of Proposition~\ref{cernov}. 
{\em Linearity}: if $f_i:\;B\to\R$, $i=1,2$, have means $P_i$ and $a,b\in\R$, then $af_1+bf_2$ has mean
$aP_1+bP_2$. If $f_1$ and $f_2$ are {\em independent}, which means that 
$\mu'(f_1^{-1}(c)\cap f_2^{-1}(d))=\mu'(f_1^{-1}(c))\cdot\mu'(f_2^{-1}(d))$ for any two values $c,d\in\R$, 
then the mean of $f_1f_2$ equals to $P_1P_2$. 
If $f$ has average at least (at most) $c$ then $f(x)\ge c$ ($f(x)\le c$) for some $x\in B$.
{\em Markov's inequality}: If $f\ge0$ has mean $P$ and $\lambda>0$, then 
$\mu'(\{x\in B\;|\;f(x)>\lambda P\})<\lambda^{-1}$. Applying it to the function $(f(x)-P)^2$ we get 
{\em \v{C}eby\v{s}ev's inequality}: If $V$, the {\em variance of $f$}, is the mean of $(f(x)-P)^2$ (where $P$ is the mean
of $f$), then 
$\mu'(\{x\in B\;|\;|f(x)-P|>\lambda\sqrt{V}\})<\lambda^{-2}$ for any $\lambda>0$.  We do not need any stronger 
result on concentration of $f$ around its mean. 

If $f:\;C\to B$ is a mapping and $\mu'$ a density on $C$, we get a density $\mu''$ on $B$ by setting 
$$
\mu''(\{x\})=\mu'(f^{-1}(x))=\sum_{c\in C,f(c)=x}\mu'(\{c\}). 
$$
We refer to this as {\em projection}. Another construction 
of more complicated densities from simpler ones takes a family of sets $B_i, i\in I$ (all of them finite), 
with a density $\mu'_1$ on $I$ and densities $\mu''_i$ on the sets $B_i$, and defines
$$
\mu'(\{(i,b)\})=\mu'_1(\{i\})\cdot\mu''_i(\{b\}),\ i\in I, b\in B_i. 
$$
Then $\mu'$ is a density on the disjoint union $\dot{\bigcup}_{i\in I}B_i$, the set of all pairs $(i,b)$ 
with $i\in I$ and $b\in B_i$. We call this construction, which generalizes to triples etc., 
{\em higher-dimensional density}. Both constructions can be combined: to define a non-uniform density
on $B$, one takes a higher-dimensional density, often patched from uniform densities, and projects it to $B$. 

Let us describe one such situation that we already encountered in the proof of Proposition~\ref{partition1} 
and will encounter again. 
Suppose that $\mu'$ is the higher-dimensional density
on $C=\dot{\bigcup}_{i\in I}B_i$ coming from the densities  $\mu'_1$ on $I$ and $\mu''_i$ on $B_i$,
$f:\;C\to B$ is a mapping that is injective for each fixed $i$ (each $B_i$ then can be regarded as a subset of $B$) 
and that $\mu''$ is the projection of $\mu'$ to $B$ via $f$. Then for each $A\sus B$, the value $\mu''(A)$ 
in fact equals to the average of the function $i\mapsto\mu''_i(B_i\cap A)$ with respect to $\mu'_1$. 

Important densities live on the sets of words $[k]^n$. The $n!$ permutations of $[n]$ act on 
the coordinates of $[k]^n$ and produce a partition $[k]^n=\bigcup_{r\in I}O_r$ into {\em orbits, or slices}, where $O_r$
consists of the words that have equal numbers of occurrences of each letter $j\in[k]$ and $I$ is the 
$\binom{n+k-1}{k-1}$-element set of $k$-tuples $r=(r_1,\ds,r_k)\in\N_0^k$, $\sum r_j=n$, recording these numbers.  
The {\em equal-slices density $\nu$} on $[k]^n$ is the unique density satisfying $\nu(\{x\})=\nu(\{y\})$ if 
$x,y\in O_r$ and $\nu(O_r)=\nu(O_s)$ for any $r,s\in I$. Explicitly, 
$$
\nu(\{x\})=\frac{1}{\binom{n+k-1}{k-1}\binom{n}{r_1,r_2,\ds,r_k}}
$$
for $x\in[k]^n$ with $r_j$ occurrences of $j$. We reserve the letter $\nu$ for the equal-slices densities and 
refer to the uniform and equal-slices densities as {\em $\mu$-density} and {\em $\nu$-density}, respectively.
If $S\sus[k]^n$ is a $d$-dimensional subspace and $A\sus[k]^n$ then $\nu_S(A)$ is defined as $\nu'(A')$ where 
$A'\sus[k]^d$ is the image of $A\cap S$ in the bijection between $S$ and $[k]^d$ and $\nu'$ is the equal-slices 
density on $[k]^d$ (this in general differs from $\nu(A\cap S)/\nu(S)$, whereas for $\mu$-density 
both ways of relativising density to subspaces give the same result). 

A slice $O_r$, $r=(r_1,\ds,r_k)$, is {\em degenerate} if $r_j=0$ for some $j$, and is {\em non-degenerate} 
else (then each letter $j\in[k]$ occurs in the words of $O_r$); there are $\binom{n-1}{k-1}$ non-degenerate slices. 
The {\em non-degenerate equal-slices density} $\tilde{\nu}$ on $[k]^n$, for $n\ge k$ (else all slices are degenerate), 
is obtained from $\nu$ by setting $\tilde{\nu}(O_r)=0$ for every degenerate slice and rescaling $\nu$ accordingly 
(by the factor $(1-\nu(D))^{-1}$ where $D$ 
is the union of degenerate slices) on the union of non-degenerate slices. So
$$
\tilde{\nu}(\{x\})=\frac{1}{\binom{n-1}{k-1}\binom{n}{r_1,r_2,\ds,r_k}}
$$
if $x\in[k]^n$ has $r_j\ge1$ occurrences of $j$ for each $j\in[k]$, and $\tilde{\nu}(\{x\})=0$ if $r_j=0$ for some $j$.

For $n,d,k\in\N$ and two words $y\in[d]^n,z\in[k]^d$, we define their {\em composition $y*z$} as the 
word $x\in[k]^n$ given by $x_i=z_{y_i}$, $i=1,2,\ds,n$. Suppose that $y\in[d]^n$ is non-degenerate (hence $n\ge d$). Clearly, 
$\{y*z\;|\;z\in[k]^d\}\sus[k]^n$ is the $d$-dimensional subspace $S(y')$, where $y'\in[k+d]^n$ is obtained from $y$ by replacing 
letter $j$ with $k+j$ (note that $S(y')$ has no fixed coordinate), and in the factorization $x=y*z$ the word $z$ is uniquely determined 
by $y$; the equality $x=y*z$ captures the way of determining $x$ by selecting a subspace $S(y')$ containing $x$ and then selecting 
`in' $S(y')$ the word $z$ corresponding to $x$. Note that if $L=L(z')\sus[k]^d$, $z'\in[k+1]^d$, is a line, then 
$\{y*z\;|\;z\in L\}=L(y*z')\sus[k]^n$ is a line too. For $n\ge d$ and $M=[d]^n\times[k]^d$, we define a new density 
$\nu'_d$ on $[k]^n$ by
$$
\nu'_d(\{x\})=\sum_{(y,z)\in M,\;y*z=x}\tilde{\nu}_1(\{y\})\cdot\tilde{\nu}_2(\{z\})
$$ 
where $\tilde{\nu}_1$ (resp. $\tilde{\nu}_2$) is the non-degenerate equal-slices density on $[d]^n$ (resp. on 
$[k]^d$); we may clearly assume that in the sum $y$ is non-degenerate. Below we show that $\nu_d'=\tilde{\nu}$. Before that 
we demonstrate that by replacing the densities $\tilde{\nu}_i$ in the definition of $\nu_d'$ with $\nu_i$ (and keeping $y$ in the sum  
non-degenerate), we obtain a density $\nu_d'$ distinct from $\nu$. Indeed, for $n=d=k=2$ and $x=11$, the two factorizations 
$11=12*11=21*11$ give $\nu_2'(\{x\})=2(3\binom{2}{1,1})^{-1}(3\binom{2}{2,0})^{-1}=1/9$, but 
$\nu(\{x\})=(3\binom{2}{2,0})^{-1}=1/3$. 

\begin{pro}\label{degslice}
Let $k,n\in\N$ and $\nu$ be the equal-slices density  and $\tilde{\nu}$ the non-degenerate equal-slices density 
on $[k]^n$.
\begin{enumerate}
\item If $m\in\N$, $j\in[k]$ and $A\sus[k]^n$ are the words with less than $m$ occurrences of $j$, then
$\nu(A)<mk/n$.
\item Let $A\sus[k]^n$, $n\ge k$, and $D\sus[k]^n$ be the union of degenerate orbits. Then (i) $\nu(D)<k^2/n$, (ii)
$\tilde{\nu}(A)=(1-\nu(D))^{-1}\nu(A)$ if $A$ consists of non-degenerate words only, and (iii) 
$|\nu(A)-\tilde{\nu}(A)|<k^2/n$ for any $A$.
\item If $n\ge d\ge k$, the above defined density $\nu'_d$ on $[k]^n$ coincides with $\tilde{\nu}$. 
\end{enumerate}
\end{pro}
\proof
1. We may set $j=k$. By the definition of $\nu$, $\nu(A)$ equals to the ratio $|M|/\binom{n+k-1}{k-1}$ where $M$ is the 
set of $k$-tuples $(r_1,\ds,r_k)\in\N_0^k$, $\sum r_i=n$, with $r_k=l<m$. Thus
$$
\nu(A)=\sum_{l=0}^{m-1}\frac{\binom{n+k-2-l}{k-2}}{\binom{n+k-1}{k-1}}=
\sum_{l=0}^{m-1}\frac{k-1}{n+k-1}\frac{\binom{n+k-2-l}{k-2}}{\binom{n+k-2}{k-2}}<\frac{mk}{n}.
$$

2. The bound on $\nu(D)$ follows from part 1 with $m=1$. The second claim is just the rescaling of $\nu$ 
defining $\tilde{\nu}$. To show the last claim set $A_1=A\cap([k]^n\backslash D)$ and $A_2=A\cap D$. 
Then (by part 1 and (ii)) $0\le\tilde{\nu}(A_1)-\nu(A_1)=\nu(D)\tilde{\nu}(A_1)\le\nu(D)<k^2/n$ and 
$0\le\nu(A_2)-\tilde{\nu}(A_2)=\nu(A_2)\le\nu(D)<k^2/n$. Since $A=A_1\cup A_2$ is a partition, 
subtraction of the two estimates gives $|\nu(A)-\tilde{\nu}(A)|<k^2/n$.

3. Let $n\ge d\ge k$, $x\in[k]^n$  be a word, $X_j\sus[n]$ for $j\in[k]$ be the positions of the 
letter $j$ in $x$ and $r_j=|X_j|$. We assume that all $r_j\ge1$ because for degenerate $x$
we clearly have $\nu'_d(\{x\})=0=\tilde{\nu}(\{x\})$ (if $y$ and $z$ are non-degenerate then so is $y*z$). The factorizations
$x=y*z$, with non-degenerate $y\in[d]^n$ and $z\in[k]^d$, 1-1 correspond to the pairs $(P,l)$ where $P$ is a partition 
of $[n]$ (a set of nonempty blocks) such that $|P|=d$ and if $B\in P$ then $B\sus X_j$ for some $j$, and 
$l:\;P\to[d]$ is a bijection. $P$ and $l$ determine $y$ and $z$ uniquely ($y_i=t\iff i\in B\in P$ with 
$l(B)=t$ and $z_i=j\iff l(B)=i$ for some $B\in P$ with $B\sus X_j$). We can generate the pairs $(P,l)$ also 
as follows. We take
all $k$-tuples $i=(i_1,\ds,i_k)\in\N^k$ with $|i|=i_1+\ds+i_k=d$, for each $i$ take all $k$ $i_j$-tuples 
$s(j)\in\N^{i_j}$, $j=1,2,\ds,k$, with $|s(j)|=r_j$, then for each $s(j)$ take all 
$\binom{r_j}{s(j)}=\frac{r_j!}{s(j)!}$ (here $s(j)!=s(j)_1!\ds s(j)_{i_j}!$) ordered partitions $(Y_{j,1},\ds,Y_{j,i_j})$
of $X_j$ with $|Y_{j,t}|=s(j)_t$, and finally we forget the orders of blocks $Y_{.,.}$ and label the $d$ blocks in each 
resulting collection in $d!$ ways 
with $1,2,\ds,d$. This way we produce each pair $(P,l)$ with multiplicity $i!=i_1!\ds i_k!$ ($i_j$ is the number of 
blocks of $P$ contained in $X_j$). Thus 
$$
\nu'_d(\{x\})=\sum_{\|i\|=k,|i|=d, \|s(j)\|=i_j,|s(j)|=r_j}\frac{d!\binom{r_1}{s(1)}\binom{r_2}{s(2)}\ds\binom{r_k}{s(k)}/i!}{\binom{n-1}{d-1}\binom{n}{s(1)s(2)\ds s(k)}
\cdot\binom{d-1}{k-1}\binom{d}{i}}
$$
where $\|\cdot\|$ is the arity of a tuple, $s(1)s(2)\ds s(k)$ means concatenation of the $i_j$-tuples into one $d$-tuple 
and the denominator gives $\tilde{\nu}_1(\{y\})\cdot\tilde{\nu}_2(\{z\})$. By cancelling the common factors in the summand we simplify the sum to
($r=(r_1,\ds,r_k)$)
$$
\frac{(n-d)!(k-1)!(d-k)!}{(n-1)!\binom{n}{r}}\sum_{\|i\|=k,|i|=d, \|s(j)\|=i_j,|s(j)|=r_j}1.
$$
The last sum equals
$$
\sum_{\|i\|=k,|i|=d}\binom{r_1-1}{i_1-1}\binom{r_2-1}{i_2-1}\ds\binom{r_k-1}{i_k-1}=\binom{n-k}{d-k}=\frac{(n-k)!}{(d-k)!(n-d)!}
$$
---we are counting $(d-k)$-element subsets $Y$ of an $(n-k)$-element set $X$ according to the sizes of 
intersections of $Y$ with blocks of a fixed partition of $X$ into blocks with sizes $r_j-1$. Hence the sum equals 
$1/\binom{n-1}{k-1}\binom{n}{r}=\tilde{\nu}(\{x\})$.
\eproof

To go from $\mu$-density to $\nu$-density, we show that if one weights $[k]^n$ on the minority of $m$ coordinates uniformly
and on the majority of remaining coordinates by equal-slices density, the resulting density is approximately $\nu$-density. 
For $k,m,n\in\N$ with $m\le n$, we define a density $\mu'_m$ on $[k]^n$ by
$$
\mu'_m(\{z\})=\sum_{M\ni(J,x,y)=z}\mu_1(\{J\})\cdot\mu_2(\{x\})\cdot\nu_1(\{y\})
$$
where $M$ consists of all triples $(J,x,y)$ with $J\sus[n]$, $|J|=m$, $x\in[k]^J$ and $y\in[k]^{[n]\backslash J}$,
a triple is projected to $[k]^n$ in the obvious way, $\mu_1$ (resp. $\mu_2=\mu_{2,J}$) is the uniform density on 
the set of $m$-element subsets of $[n]$ (resp. on $[k]^J$) and $\nu_1=\nu_{1,J}$ is the equal-slices density on 
$[k]^{[n]\backslash J}$. 

\begin{pro}\label{jejichL66}
Let $k,m,n\in\N$, $m\le n$ and $\mu'_m$ be the above defined density on $[k]^n$. Then for every set $A\sus[k]^n$,  
$$
|\mu'_m(A)-\nu(A)|\le km/n.
$$
\end{pro}
\proof
We prove the inequality, in fact a stronger one, first for $m=1$. Let $z\in[k]^n$ and $r_j$ be the number of occurrences of the letter $j\in[k]$ in $z$. By the definition of $\mu'_m$ and $\nu$,
\begin{eqnarray*}
\frac{\mu'_1(\{z\})}{\nu(\{z\})}&=&\binom{n+k-1}{k-1}\binom{n}{r_1,r_2,\ds,r_k}
\sum_{j=1, r_j\ge1}^k\frac{r_j/kn}{\binom{n+k-2}{k-1}\binom{n-1}{r_1,\ds,r_j-1,\ds,r_k}}\\
&=&\frac{\binom{n+k-1}{k-1}}{\binom{n+k-2}{k-1}}\sum_{j=1, r_j\ge1}^k\frac{1}{k}
\le\frac{\binom{n+k-1}{k-1}}{\binom{n+k-2}{k-1}}=1+\frac{k-1}{n}.
\end{eqnarray*}
So $|\mu'_1(\{z\})-\nu(\{z\})|\le\frac{k-1}{n}\nu(\{z\})$. Summing over $z\in A$ and using triangle inequality 
we deduce that 
$${\textstyle
|\mu'_1(A)-\nu(A)|\le\frac{k-1}{n}\nu(A)<\frac{k}{n}.
}
$$

We derive from this that $|\mu'_m(A)-\mu'_{m-1}(A)|\le k/n$ for every $A\sus[k]^n$ and $m\ge2$. 
The inequality $|\mu'_m(A)-\nu(A)|\le mk/n$ then follows by induction and triangle inequality. 
Let $m\ge2$ and $A\sus[k]^n$. We partition the set of triples $M$ defining $\mu'_{m-1}$ by the equivalence 
$(J,x,y)\sim(J',x',y')$ iff $J=J'$ and $x=x'$. So (projecting $(J,x,y)$ to $[k]^n$ when needed)
$$
\mu'_{m-1}(A)=\sum_{B\in M/\!\sim}\mu_1(\{J\})\cdot\mu_2(\{x\})\sum_{B\ni(J,x,y)\in A}\nu_1(\{y\}).
$$
We replace in each inner sum the equal-slices density $\nu_1$ on 
$[k]^{[n]\backslash J}$ (now $|J|=m-1$) with the density $\mu'_1$, and this changes the total sum to, say, $\mu''(A)$. 
Summing the changes of inner sums, the result for $m=1$ gives that $|\mu''(A)-\mu'_{m-1}(A)|\le(k-1)/(n-m+1)$. A moment
of reflection reveals that the change of $\nu_1$ on each $[k]^{[n]\backslash J}$ to $\mu'_1$ gives an equivalent, 
only a more complicated, way of counting $\mu'_m(A)$ (it boils down to the identity 
$\binom{n}{m}=\frac{n-m+1}{m}\binom{n}{m-1}$). Thus $\mu''(A)=\mu'_m(A)$ and 
$|\mu'_m(A)-\mu'_{m-1}(A)|\le (k-1)/(n-m+1)\le k/n$ as needed, because we may assume that $n\ge km$ (else the result 
holds trivially). 
\eproof

Propositions~\ref{degslice} and \ref{jejichL66} show that we may replace $\mu$-density in Theorem~\ref{denshalejewe} 
with $\tilde{\nu}$-density. Using this we derive Proposition~\ref{pravddhjthm}, the key strengthening of 
Theorem~\ref{denshalejewe}. 

For $k,m,n\in\N$, $J\sus[n]$ with $|J|=m$ and $y\in[k]^{[n]\backslash J}$, we denote by $S_{J,y}$
the $m$-dimensional subspace of $[k]^n$ that has $J$ as the set of free coordinates and elsewhere is determined by $y$:  
$x\in S_{J,y}\iff x_i=y_i$ for every $i\in[n]\backslash J$.

\begin{pro}\label{eqsldhjthm}
Let $k\in\N$, $k\ge2$, be given and assume Theorem~\ref{denshalejewe} for $k$. It follows that for every $\de>0$ 
there is an $n_0\in\N$ such that for every $n>n_0$ every set $A\sus[k]^n$ with $\tilde{\nu}(A)\ge\de$ contains a line.
\end{pro}
\proof
Let $\de$ be given. We take the $n_0$ of Theorem~\ref{denshalejewe} corresponding to uniform density 
$\de/3$ and set $m=n_0+k$. Suppose that $n>3km/\de=3k(n_0+k)/\de$ and that $A\sus[k]^n$ has $\tilde{\nu}(A)\ge\de$. 
By part 2 (iii) of Proposition~\ref{degslice}, $\nu(A)\ge2\de/3$. By Proposition~\ref{jejichL66} and the definition of 
density $\mu_m'$ before it, there exists an $m$-dimensional 
subspace $S=S_{J,y}$ of $[k]^n$, $J\sus[n]$ with $|J|=m$, such that $\mu_S(A)\ge\nu(A)-km/n\ge\de/3$. By the choice of 
$m$ and Theorem~\ref{denshalejewe}, there is a line in $[k]^n$ contained in $A\cap S$. 
\eproof

\begin{pro}\label{pravddhjthm}
Let $k\in\N$, $k\ge2$, be given and assume Theorem~\ref{denshalejewe} for $k$. It follows that for every 
$\de>0$ there exist an $n_0\in\N$ and a $\theta>0$ such that if $n>n_0$ and 
$A\sus[k]^n$ has $\nu(A)\ge\de$, then the set $M\sus[k+1]^n$ of lines contained in $A$ has 
$$
\nu(M)\ge\theta.
$$
\end{pro}
\proof
Let $\de$ be given. We take the $n_0$ of Proposition~\ref{eqsldhjthm} corresponding to the $\tilde{\nu}$-density $\de/2$ and set $d=n_0+1$. Suppose that $n>d+4k^2/\de$ and $A\sus[k]^n$ has $\nu(A)\ge\de$.
By part 2 (iii) of Proposition~\ref{degslice}, $\tilde{\nu}(A)\ge3\de/4$. For $y\in[d]^n$ we define 
$C_y=\{z\in[k]^d\;|\;y*z\in A\}$ (recall the composition of words $*$ introduced before Proposition~\ref{degslice}).
Let $B\sus[d]^n$ be the set of (non-degenerate) words $y$ such that $\tilde{\nu}_2(C_y)\ge\de/2$. By part 3 of 
Proposition~\ref{degslice} (applied with $k$), $\tilde{\nu}(A)\ge3\de/4$ and the definition of $B$ imply that 
$\tilde{\nu}_1(B)\ge\de/4$. 
Deleting degenerate words (they are irrelevant for $\tilde{\nu}_2$ anyway), we may assume that all words in 
every $C_y$ are non-degenerate. Consider the set
$$
M=\{x'\in[k+1]^n\;|\;x'=y*z',y\in B,z'\in[k+1]^d,L(z')\sus C_y\}.
$$
These are lines contained in $A$: $L(x')\sus A$ for every $x'\in M$. By the choice of $d$, for each $y\in B$
the set $C_y\sus[k]^d$ contains a line $L(z')$ and (due to the purge on $C_y$) $z'$ is non-degenerate. 
Let $\tilde{\nu}_2'$ (resp. $\tilde{\nu}'$) be the non-degenerate equal-slices density on $[k+1]^d$ (resp. on $[k+1]^n$). 
Since $\tilde{\nu}_1(B)\ge\de/4$ 
and for each $y\in B$ there is at least one non-degenerate $z'\in[k+1]^d$ with $y*z'\in M$, giving contribution 
at least $\tilde{\nu}_2'(\{z'\})>d^{-k}(k+1)^{-d}$, by part 3 of Proposition~\ref{degslice} (applied with $k+1$) 
we see that
$$
\tilde{\nu}'(M)>\frac{\de}{4d^k(k+1)^d}.
$$
By part 2 ((i) and (ii)) of  Proposition~\ref{degslice} (and since $n>4k^2$ and $k\ge2$), the desired 
lower bound $\nu(M)>(1-(k+1)^2/n)\tilde{\nu}'(M)>(\de/9)d^{-k}(k+1)^{-d}=\theta>0$ follows.  
\eproof

To go from $\nu$-density to $\mu$-density, we show that if one weights $[k]^n$ on the minority of $m$ coordinates 
by $\nu$-density and on the majority of remaining coordinates uniformly, the resulting density is approximately 
$\mu$-density. We prove it in greater generality with any density $\mu'$ on the minority of $m$ coordinates. 
For $k,m,n\in\N$ with $m\le n$ and a density $\mu'$ on $[k]^m$, we define a density $\mu'_m$ on $[k]^n$ by
$$
\mu'_m(\{z\})=\sum_{M\ni(\si,x,y)=z}\mu_1(\{\si\})\cdot\mu'(\{x\})\cdot\mu_2(\{y\})
$$
where $M$ consists of all triples $(\si,x,y)$ with $\si:\;[m]\to[n]$ an injection, $x\in[k]^m$ and 
$y\in[k]^{[n]\backslash\si([m])}$, a triple $(\si,x,y)$ is projected to $[k]^n$ by setting $z_{\si(i)}=x_i$ for $i\in[m]$ and 
$z_i=y_i$ for $i\in[n]\backslash\si([m])$, $\mu_1$ (resp. $\mu_2=\mu_{2,\si}$) is the uniform density on the set of injections  
from $[m]$ to $[n]$ (resp. on $[k]^{[n]\backslash\si([m])}$) and $\mu'$ is the given density on $[k]^m$.

\begin{pro}\label{prumdens}
Let $k,m,n\in\N$ and $\eta>0$ be such that $m\le n^{1/4}$ and $n\ge(12k/\eta)^{12}$, $\mu'$ be a density on
$[k]^m$ and $\mu'_m$ be the above corresponding density on $[k]^n$. Then for every set $A\sus[k]^n$, 
$$
|\mu_m'(A)-\mu(A)|<\eta.
$$
\end{pro}
\proof
It suffices to consider only $\mu'=\mu'_u$ given, for some $u\in[k]^m$, by $\mu'(\{u\})=1$ and $\mu'(\{x\})=0$ for
$x\ne u$, because any density $\mu'$ on $[k]^m$ is a convex combination of these densities, 
$\mu'=\sum_u\lambda_u\mu'_u$ ($\lambda_u\ge0$, $\sum_u\lambda_u=1$), and $\mu'_m=\sum_u\lambda_u\mu_{u,m}'$; 
the general result follows by the triangle inequality. 

We fix words $u\in[k]^m$ and $z\in[k]^n$ such that $z$ has between $n/k-n^{2/3}$ and $n/k+n^{2/3}$ occurrences 
of each letter $j\in[k]$ (by Proposition~\ref{cernov}, only very few $z$ are not like this). 
If $p$ (resp. $q$) is the minimum (resp. maximum) number of occurrences 
of a letter $j$ in $z$ (clearly $p>m$) then 
$$
\left(\frac{p-m}{n}\right)^mk^{m-n}\le\mu_m'(\{z\})\le\left(\frac{q}{n-m}\right)^mk^{m-n}
$$
because $\mu_1(\{\sigma\})=1/n(n-1)\ds(n-m+1)$ lies between $n^{-m}$ and $(n-m)^{-m}$, 
the number of $\sigma$ satisfying $u_i=z_{\sigma(i)}$ for $i\in[m]$ is at least $(p-m+1)^m$ and at most $q^m$ 
($\sigma$ determines $x$ and $y$) and $\mu_2(\{y\})=k^{m-n}$. Since $\mu(\{z\})=k^{-n}$, $n/k-n^{2/3}\le p\le q\le n/k+n^{2/3}$ and $m\le n^{1/4}$, we have
$$
(1-2kn^{-1/3})^m<\frac{\mu_m'(\{z\})}{\mu(\{z\})}\le\left(\frac{1+kn^{-1/3}}{1-n^{-3/4}}\right)^m<(1+2kn^{-1/3})^m.
$$ 
Since $1-\de<e^{-\de}<1-\de/2$ and $1+\de<e^{\de}<1+2\de$ if $\de\in(0,\frac{1}{2})$, we deduce that
$\frac{\mu_m'(\{z\})}{\mu(\{z\})}$ lies in $(1-4kmn^{-1/3},1+4kmn^{-1/3})$ provided that 
$4kmn^{-1/3}<\frac{1}{2}$. This is true as $4kmn^{-1/3}\le 4kn^{-1/12}\le\eta/3<\frac{1}{2}$. So 
$$
\left|\mu_m'(\{z\})/\mu(\{z\})-1\right|<\eta/3.
$$

Let $[k]^n=B\cup C$, where $B$ are the words meeting the condition on occurrences of letters 
and $C$ are the remaining words. By Proposition~\ref{cernov}, $\mu(C)<kn^{-1/3}<\eta/3$. Since 
$\mu_m'(\{z\})>(1-\eta/3)\mu(\{z\})$ for $z\in B$, we have 
$\mu_m'(B)>(1-\eta/3)\mu(B)>(1-\eta/3)^2>1-2\eta/3$ and $\mu_m'(C)<2\eta/3$. We conclude that
\begin{eqnarray*}
|\mu_m'(A)-\mu(A)|&\le&\sum_{z\in A\cap B}|\mu_m'(\{z\})-\mu(\{z\})|+|\mu_m'(A\cap C)-\mu(A\cap C)|\\
&<&\sum_{z\in A\cap B}\mu(\{z\})(\eta/3)+2\eta/3\\
&\le&\eta.
\end{eqnarray*} 
\eproof

\noindent
We apply Proposition~\ref{prumdens} to three densities $\mu'$ on $[k]^m$, all invariant to permuting the 
$m$ coordinates. The definition of $\mu_m'$ then simplifies, as one can put the $\sigma$ with the common $m$-element image 
$J\sus[n]$ together and sum over the triples $(J,x,y)$. The first application with   
$\mu'(A)=\mu_B(A)$, where $A\sus[k]^m$ and $B=[k-1]^m$, gives Proposition~\ref{mucarka}.
In the other two applications of Proposition~\ref{prumdens}, $\mu'$ is the 
equal-slices density on $[k]^m$, respectively the density given by $\mu'(A)=\nu'(A\cap[k-1]^m)$ where $\nu'$ 
is the equal-slices density on $[k-1]^m$, and we get the next proposition, for which we introduce the following notation. 
The {\em truncation} $S'\sus S\sus[k]^n$ of an $m$-dimensional subspace $S$ is obtained by forbidding $k$ as the value 
of $x\in S$ on the free coordinates; $S'$ 1-1  corresponds 
with $[k-1]^m$. For $A\sus[k]^n$ we define $\nu_{S'}(A)$ as $\nu'(A')$ where $A'$ is the image of $A\cap S'$ 
in the bijection between $S'$ and $[k-1]^m$ and $\nu'$ is the equal-slices density on $[k-1]^m$. 

\begin{pro}\label{jejichC64i65}
Let $\de,\eta>0$ and $k,m,n\in\N$ satisfy $m\le n^{1/4},n\ge(12k/\eta)^{12}$ and $A\sus[k]^n$ be a set
with $\mu(A)=\de$. Then the (uniform) averages of the functions $S\mapsto\nu_S(A)$ and $S\mapsto\nu_{S'}(A)$, over all 
subspaces $S=S_{J,y}$, $J\sus[n]$, $|J|=m$, and words $y\in[k]^{[n]\backslash J}$, are both at least $\de-\eta$.
\end{pro}

\noindent
To deduce Proposition~\ref{inssens}, we need Propositions~\ref{degslice}, \ref{jejichL66}, \ref{pravddhjthm} and 
\ref{jejichC64i65}.

\subsection{Derivation of Proposition~\ref{inssens}}

In this subsection we fix a $k\in\N$, $k\ge3$, assume that Theorem~\ref{denshalejewe} holds for $k-1$ (and every 
$\de>0$) and deduce from this Proposition~\ref{inssens} for $k$. We proceed in three steps. First we show that for any 
positively $\mu$-dense $A\sus[k]^n$ there is a subspace $S\sus[k]^n$ such that either $A$ gets on $S$ $\nu$-denser (which 
gives the desired density increment at once), or $A$ gets positively $\nu$-dense on the truncation $S'$ of $S$ 
(recall that $S'$ has forbidden $k$ on the free coordinates) while losing not too much $\nu$-density on the whole $S$. 
In the crucial second step we obtain, 
assuming the second alternative and that $A$ is free of lines, a $\nu$-density increment of $A$ on a $k$-set $D$, 
an increment large enough to make up for the previous loss. In the third step we convert the $\nu$-density increment of $A$ on $D$ to 
a $\mu$-density increment.

\begin{pro}\label{1ststep}
Let $\de,\eta>0$ and $m,n\in\N$ satisfy $\eta\le\de/4$, $m\le n^{1/4}$, $n\ge(12k/\eta)^{12}$.
Then for every set $A\sus[k]^n$ with $\mu(A)=\de$ there exists an $m$-dimensional subspace
$S\sus[k]^n$ such that 1 or 2 holds:
\begin{enumerate}
\item $\nu_S(A)\ge\de+\eta=\mu(A)+\eta$;
\item $\nu_S(A)\ge\de-4\eta\de^{-1}=\mu(A)-4\eta\de^{-1}$ and $\nu_{S'}(A)\ge\de/4$, where $S'\sus S$ 
is the truncation of $S$ with values on the free coordinates lying in $[k-1]$.
\end{enumerate}
\end{pro}
\proof
We take uniformly the subspaces $S=S_{J,y}$, as described in Proposition~\ref{jejichC64i65}. Let $M$ be the set 
of $S$ with $\nu_S(A)<\de-4\eta/\de$ and $N$ be the set of $S$ with $\nu_{S'}(A)<\de/4$. We assume that 1 does not
hold, so $\nu_S(A)<\de+\eta$ for every $S$, and show that then 2 holds. If $\mu(M)\ge\de/2$ then the average of 
$\nu_S(A)$ over $S$ is at most
$$
(1-\de/2)(\de+\eta)+(\de/2)(\de-4\eta/\de)=\de+(1-\de/2)\eta-2\eta<\de-\eta,
$$
contradicting Proposition~\ref{jejichC64i65}. So $\mu(M)<\de/2$. Similarly, if $\mu(N)\ge1-\de/2$ then the average 
of $\nu_{S'}(A)$ over $S$ is at most
$$
\de/2+(1-\de/2)(\de/4)<3\de/4\le\de-\eta,
$$
again contradicting Proposition~\ref{jejichC64i65}. So $\mu(N)<1-\de/2$. Hence there is a subspace $S=S_{J,y}$ not in 
$M\cup N$ and 2 holds.
\eproof

For $x\in[k]^m$ (we have replaced $n$ by $m$ to indicate that we move into $S$) and $j\in[k-1]$, we denote, as before, 
by $x(j)$ the word obtained from $x$ by changing all $k$s to $j$s. For a set $A_1\sus[k]^m$ and $j\in[k-1]$, we define 
$$
C_j=\{x\in[k]^m\;|\;x(j)\in A_1\}\ \mbox{ and }\ C=\bigcap_{j=1}^{k-1}C_j\sus[k]^m.
$$
Note that each $C_j$ is $jk$-insensitive and that, crucially, if $A_1$ contains no line then $A_1\cap C\sus[k-1]^m$. Indeed, 
if $x\in A_1\cap C$ had an occurrence of $k$, then $\{x\}\cup\{x(j)\;|\;j\in[k-1]\}$ would be a line in $[k]^m$ 
contained in $A_1$.

\begin{pro}\label{heart}
For every $\de_1>0$ there is an $m_0\in\N$ and a $\theta>0$
such that the following holds. If $m>m_0$ and $A_1\sus[k]^m$ contains no line, $\nu(A_1)\ge\de_1$ and 
(measured in $[k-1]^m$) $\nu(A_1\cap[k-1]^m)\ge\de_1/4$, then there is a $k$-set $D\sus[k]^m$ satisfying
$$
\nu(A_1\cap D)\ge\nu(A_1)\nu(D)+\de_1\theta/2k\ge\de_1\nu(D)+\de_1\theta/2k.
$$
\end{pro}
\proof 
Let $\de_1$ be given. Applying Theorem~\ref{denshalejewe} for $k-1$ and Proposition~\ref{pravddhjthm}, we take an $m_0$ and a 
$\theta>0$ such that 
if $m>m_0$ then for every set $B\sus[k-1]^m$ with $\nu(B)\ge\de_1/4$ the set $M\sus[k]^m\backslash[k-1]^m$ 
of lines contained in $B$ has $\nu(M)\ge\theta$; we also assume $m_0$ so big that $m>m_0$ implies $k/\theta m<\de_1/2$.
Now let $m>m_0$ and $A_1\sus[k]^m$ be as stated, with the above defined sets $C_j$ and $C$; for convenience we denote
$\de_1=\nu(A_1)$. By the assumptions we may take as $B$ the set $B=A_1\cap[k-1]^m$. 
The lines $M$ in $B$ then 1-1 correspond to the words in $C\backslash[k-1]^m$. Hence $\nu(C\backslash[k-1]^m)\ge\theta$. 
We observed above that $C\backslash[k-1]^m$ is disjoint to $A_1$. Therefore using part 1 of Proposition~\ref{degslice} we get 
$$
\nu(A_1\cap C)\le k/m<\theta\de_1/2\le(\de_1/2)\nu(C). 
$$
For $j\in[k]$ we set $D^{(j)}=C_1\cap\ds\cap C_{j-1}\cap([k]^m\backslash C_j)$; $D^{(1)}=[k]^m\backslash C_1$ and
$D^{(k)}=C$. Thus $[k]^m=\bigcup_{j=1}^kD^{(j)}$ is a partition. By $\nu(A_1)=\de_1$ and 
$\nu(A_1\cap D^{(k)})\le(\de_1/2)\nu(D^{(k)})$,
\begin{eqnarray*}
\nu(A_1\cap(D^{(1)}\cup\ds\cup D^{(k-1)}))&\ge&\de_1-(\de_1/2)\nu(D^{(k)})\\
&=&\de_1(1-\nu(D^{(k)}))+(\de_1/2)\nu(D^{(k)})\\
&\ge&\de_1\nu(D^{(1)}\cup\ds\cup D^{(k-1)})+\de_1\theta/2.
\end{eqnarray*}
Thus $\nu(A_1\cap D^{(j)})\ge\de_1\nu(D^{(j)})+\de_1\theta/2(k-1)$ for some $j\in[k-1]$. We set, for this $j$, 
$D_i=C_i$ for $i<j$, $D_j=[k]^m\backslash C_j$ and $D_i=[k]^m$ for $i>j$. Clearly, each $D_i$ is $ik$-insensitive. 
The $k$-set $D=\bigcap_{i=1}^{k-1}D_i=D^{(j)}$ satisfies the displayed inequality.
\eproof

\noindent
This is the heart of the proof of Theorem~\ref{denshalejewe}, transmuting the inductive assumption 
on the level $k-1$ in a density increment on the level $k$. The quantities $m_0=m_0(\de_1)$ and $\theta=\theta(\de_1)$
come from the validity of Theorem~\ref{denshalejewe} for $k-1$. In particular, note that $\theta$ can be 
assumed nondecreasing in $\de_1$ (it is obvious from the proof but perhaps is not so clear from the statement).

\begin{pro}\label{3rdstep}
Let $\be,\de_2\in(0,1)$, $m,r\in\N$, $\be\ge kr/m$ and let $A_2\sus D\sus[k]^m$ be sets satisfying 
$\nu(A_2)\ge\de_2\nu(D)+3\be$. Then there exists a subspace $V\sus[k]^m$ with dimension $r$ such that
$$
\mu_V(A_2)\ge\de_2\mu_V(D)+\be.
$$
\end{pro}
\proof
The average of $\mu_V(A_2)-\de_2\mu_V(D)$ over all subspaces $V=S_{J,y}$, with $J\sus[m]$, $|J|=r$, taken uniformly and 
$y\in[k]^{[m]\backslash J}$ taken according to $\nu$-density, equals $\mu_r'(A_2)-\de_2\mu_r'(D)$ where 
$\mu'_r$ is the density on $[k]^m$ introduced before Proposition~\ref{jejichL66}. By Proposition~\ref{jejichL66} 
and the assumptions this is at least $(\nu(A_2)-\be)-\de_2(\nu(D)+\be)=\nu(A_2)-\de_2\nu(D)-2\be\ge\be$. 
Thus a subspace $V=S_{J,y}$ exists that satisfies the displayed inequality.
\eproof

We prove Proposition~\ref{inssens}. Let $r\in\N$ and $\de>0$ be given ($k\ge3$ is fixed) and suppose that 
$A\sus[k]^n$ contains no line, $\mu(A)\ge\de$ and $n>n_0$; we specify a bound on $n_0$ at the end.
We set $\de_1=\de/2$ and take the $m_0=m_0(\de_1)$ and $\theta=\theta(\de_1)$ of Proposition~\ref{heart}.
Let $\eta=\de^2\theta/32k$ and $m=\lfloor n^{1/4}\rfloor$. Note that $\eta<\de/4$, $\de-4\eta\de^{-1}>\de_1$, 
$m\le n^{1/4}$ and $\de/4>\de_1/4$. By Proposition~\ref{1ststep}, applied for $\de=\mu(A),\eta$, $m$ and $A$, 
if $n\ge(12k/\eta)^{12}$ then there is an $m$-dimensional subspace $S\sus[k]^n$ satisfying alternative 1 or alternative 2. 
We denote by $A_1\sus[k]^m$ the image of $A\cap S$ in the bijection between $S$ and $[k]^m$. We first consider alternative 1. 
So $\nu(A_1)\ge\mu(A)+\eta$. 
For $n$ large enough so that $\eta m/3k\ge r$, Proposition~\ref{3rdstep}, applied for $\be=\eta/3$, $\de_2=\mu(A)$, 
$m,r$, $A_2=A_1$ and $D=[k]^m$, provides an $r$-dimensional subspace $V\sus[k]^m$ on which 
$\mu_V(A_1)\ge\mu(A)+\eta/3$. We achieved a $\mu$-density increment of $A$ on the $k$-set $D=W$ in the 
$r$-dimensional subspace $W\sus[k]^n$ that is the image of $V$ in the bijection between $[k]^m$ and $S$, with the 
increment $\ga=\eta/3$. Clearly, $\mu_W(D)=1>\ga$.

Let $S\sus[k]^n$ satisfy alternative 2 of Proposition~\ref{1ststep}. 
So $\nu(A_1)\ge\mu(A)-4\eta/\mu(A)>\de_1$ and $\nu(A_1\cap[k-1]^m)>\de/4>\de_1/4$. By Proposition~\ref{heart}, 
applied for $\de_1$ and $A_1$, for large enough $n$ (so that $m>m_0$) there is a $k$-set $D_1\sus[k]^m$ for which 
\begin{eqnarray*}
\nu(A_1\cap D_1)&\ge&\nu(A_1)\nu(D_1)+\de_1\theta/2k\ge\mu(A)\nu(D_1)-4\eta\de^{-1}+\de_1\theta/2k\\
&>&\mu(A)\nu(D_1)+\de\theta/8k.
\end{eqnarray*}
We apply Proposition~\ref{3rdstep} with $\be=\de\theta/24k$, $\de_2=\mu(A)$, 
$m,r$, $A_2=A_1\cap D_1$ and $D_1$. For large enough $n$ (so that $\be\ge kr/m$) it provides an $r$-dimensional subspace 
$V\sus[k]^m$ with $\mu_V(A_2)\ge\mu(A)\mu_V(D_1)+\be$. Note that $D_1\cap V$ is a $k$-set in $V$.
We achieved a $\mu$-density increment of $A$ on the $k$-set $D=c(D_1\cap V)$ in the $r$-dimensional subspace 
$W=c(V)\sus[k]^n$, where $c$ is the bijection between $[k]^m$ and $S$, with the 
increment $\ga=\be=\de\theta/24k$. Clearly, $\mu_W(D)\ge\be=\ga$.

To summarize and integrate both cases, we see that for given $r\in\N$, $\de\in(0,1)$ and any $n>n_0$, for any 
set $A\sus[k]^n$ containing no line and with $\mu(A)\ge\de$ there is an $r$-dimensional subspace $W\sus[k]^n$ 
and a $k$-set $D\sus W$ in $W$ such that 
$\mu_W(D)\ge\ga$ and $\mu_W(A\cap D)\ge\mu(A)\mu_W(D)+\ga$ (hence $\mu_D(A)\ge\mu(A)+\ga$), with the desired density increment 
$\ga=\min(\eta/3,\be)=\eta/3=\de^2\theta/96k$. We observed above that $\theta$ is nondecreasing in $\de_1=\de/2$ 
and so $\ga$ is nondecreasing in $\de$. Finally, the argument shows that the sufficient $n_0$ to take is, 
for $\eta=\de^2\theta/32k$,
$$
n_0=\lceil(12k/\eta)^{12}+(3kr/\eta)^4+m_0^4+(24k^2r/\de\theta)^4\rceil
$$
where $m_0=m_0(\de/2)$ and $\theta=\theta(\de/2)$ are the quantities of Proposition~\ref{heart}, 
guaranteed by Theorem~\ref{denshalejewe} for $k-1$. This concludes the derivation of Proposition~\ref{inssens}. 

The proof of Theorem~\ref{denshalejewe}, the density Hales--Jewett theorem, and consequently of 
Theorem~\ref{szem}, Szemer\'edi's theorem, is complete. 

\section{Concluding remarks and thoughts}

In writing this text we were motivated also by the last sentence of the abstract in \cite{poly}: 
``Our proof is surprisingly simple: indeed, it gives arguably the simplest known proof of Szemer\'edi's theorem.''
How simple/long is then Polymath's proof of Szemer\'edi's theorem? The article \cite{poly} has
44 pages but the proof of Theorem~\ref{denshalejewe} only starts after 32 pages in Section 7 and takes about
8 pages, during which it draws on various results and concepts obtained in the preceding part. The original article of 
Szemer\'edi \cite{szem} has 46 pages and Furstenberg's ergodic paper \cite{furs} 52. In the book of Moreno and 
Wagstaff, Jr. \cite[Chapter 7]{more_wags}, one of the few (if not the only one) monographs or textbooks presenting Szemer\'edi's
combinatorial proof of his theorem, the proof takes 38 pages, and in the write-up of Tao \cite{tao_writeup}
about 26. An article of Tao \cite{tao_szem} of 49 pages gives a proof of Szemer\'edi's theorem based on a combination of 
ergodic methods and the approach of Gowers \cite{gowe_01}. Towsner \cite{tows} gives a (not quite self-contained) 
model-theoretic proof of Szemer\'edi's theorem on 10 pages. (This list of proofs of Theorem~\ref{szem} in the literature or
on the Internet is far from exhaustive.) Our present write-up, a reshuffled and pruned form of Polymath's proof 
\cite{poly}, demonstrates that it is possible to write down a self-contained combinatorial proof of Szemer\'edi's theorem
well under 20 pages, which justifies the quoted sentence. Of course, it is even a proof of a stronger theorem, the density 
Hales--Jewett theorem. 

As for the correctness of the proof in \cite{poly}, we pointed in the remark after the proof of Proposition~\ref{partition1} 
a probably overlooked lower bound factor in \cite[Lemma 8.1]{poly}, but this is trivial to repair (which we did) and we 
did not notice in \cite{poly} anything more serious than that. In recent years formal proofs of various popular theorems
were worked out, for example, for the Prime Number Theorem (Avigad et al. \cite{avig_spol}, Harrison \cite{harr1}), Dirichlet's theorem 
on primes in arithmetic progression (Harrison \cite{harr2}) or Jordan's curve theorem (Hales \cite{hale}). Szemer\'edi's theorem
is known for logical intricacy of its proof---an interesting project in formal proofs may be to produce a formal version for it or, 
for this matter, for the proof of the density Hales--Jewett theorem. 

Many arguments of the proof in \cite{poly} as we present them are simple instances of the probabilistic method reasoning
(see Alon and Spencer \cite{alon_spen}), but we evade words `probability', `random' or `randomly' in our write-up
(in \cite{poly} the last two words appear more than 90 times). We prefer the terminology of densities instead, to emphasize that we give 
in all cases explicit definitions and constructions of the densities (i.e., probability measures) used, which is not quite done in \cite{poly}. 
We consider it important, for the sake of rigorousness of the whole approach, to give these explicit definitions. 
For illustration consider the identity in part 3 of Proposition~\ref{degslice}, for which 
we gave a verificational proof. The original proof of Polymath \cite[pp. 1297--8]{poly}, more elegant, is free of calculations and is 
based (in our terminology) on representing the non-degenerate equal-slices density $\tilde{\nu}$ on $[k]^n$, 
$k\le n$, as a projection of a higher-dimensional density built from uniform densities. Informally (\cite[p. 1295]{poly}): 
a $\tilde{\nu}$-random word $x$ arises by selecting  $n$ points $q_1,\ds,q_n$ around a circle in a random order, putting randomly 
$k$ delimiters $r_1,\ds,r_k$ in some $k$ distinct gaps out of the $n$ gaps determined by the $n$ points $q_i$, and 
then reading the positions of the letter $j\in[k]$ in $x$ in the indices $i$ of the points $q_i$ lying between $r_j$ and the delimiter clockwisely preceding $r_j$. Formally (not an exact translation): for $x\in[k]^n$, let
$$
\mu'(\{x\})=\sum_{(\pi,B')=x}\mu_1(\{\pi\})\cdot\mu_2(\{B'\})
$$
where $\pi$ run through the $n!$ permutations of $[n]$, $B'$ run through the $k\binom{n}{k}$ pointed $k$-element subsets 
of $[n]$ (the pairs $(B,b)$ with $b\in B\sus[n]$, $|B|=k$), $\mu_i$ are uniform densities and $(\pi,B')$ projects
to $[k]^n$ as follows. If $\pi=a_1a_2\ds a_n$ and $B'=(B,b)$ with $B=\{b_1<b_2<\ds<b_k\}$ and $b=b_t$, 
we project $(\pi,B')$ to $x\in[k]^n$ by setting, for $j=1,2,\ds,k$, $x_{a_i}=j$ exactly for the terms $a_i$ in $\pi$ with $i$ 
in the interval $b_{t+j-1}\le i<b_{t+j}$, where the indices are taken modulo $k$ and the interval 
$b_k\le i<b_{k+1}=b_1$ is $[b_k,n]\cup[1,b_1)$. It is immediate to show that $\mu'=\tilde{\nu}$. 

In conclusion, we want to remark that the use of non-uniform densities $\nu$ and $\tilde{\nu}$ on words and their interplay with 
the uniform density is a really interesting and combinatorially beautiful feature of Polymath's proof \cite{poly}.  



\begin{thebibliography}{10}
\bibitem{ajta_szem} M. Ajtai and E. Szemer\'edi, Sets of lattice points that form no squares, 
{\em Studia Sci. Math. Hung.} {\bf 9} (1974), 9--11 (1975).
\bibitem{alon_spen} N. Alon and J.\,H. Spencer, {\em The Probabilistic Method}, John Wiley $\&$ Sons,
2008 (3rd edition). 
\bibitem{avig_spol} J. Avigad, K. Donnelly, D. Gray and P. Raff, A formally verified proof of the prime number theorem, 
{\em ACM Trans. Comput. Log.} {\bf 9} (2008) 23 pp. 
\bibitem{bara_soly} I. B\'ar\'any and J. Solymosi (Editors), {\em An Irregular Mind. Szemer\'edi is 70}, 
Springer, Berlin, 2010.
\bibitem{furs} H. Furstenberg, Ergodic behavior of diagonal measures and a theorem of Szemer\'edi
on arithmetic progressions, {\em J. Analyse Math.} {\bf 31} (1977), 204--256.
\bibitem{furs_katz_89} H. Furstenberg and Y. Katznelson, A density version of the Hales--Jewett theorem
for $k=3$, {\em Discrete Math.} {\bf 75} (1989), 227--241.
\bibitem{furs_katz} H. Furstenberg and Y. Katznelson, A density version of the Hales--Jewett theorem, 
{\em J. Analyse Math.} {\bf 57} (1991), 64--119.
\bibitem{gowe_01} W.\,T. Gowers, A new proof of Szemer\'edi's theorem, {\em Geom. Funct. Anal.} {\bf 11} (2001), 
465--588. 
\bibitem{gowe} W.\,T. Gowers, Polymath and the density Hales--Jewett theorem, pp. 659--688 in \cite{bara_soly}.
\bibitem{hale_jewe} A.\,W. Hales and R.\,I. Jewett, Regularity and positional games, {\em Trans. AMS} 
{\bf 106(2)} (1963), 222--229.
\bibitem{hale} T.\,C. Hales, The Jordan curve theorem, formally and informally, {\em Amer. Math. Monthly} 
{\bf 114} (2007), 882--894.
\bibitem{harr1} J. Harrison, Formalizing an analytic proof of the prime number theorem, {\em J. Automat. Reason.}
{\bf 43} (2009), 243--261. 
\bibitem{harr2} J. Harrison, A formalized proof of Dirichlet's theorem on primes in arithmetic progression, {\em J. Formaliz. 
Reason.} {\bf 2} (2009), 63--83. 
\bibitem{klaz_I} M. Klazar, Diophantine equation $ax^n-by^n=c$. I, {\em KAM--DIMATIA Series}, 2010-971, 39 pp. 
\bibitem{klaz_II} M. Klazar, Analytic and combinatorial number theory I, {\em KAM--DIMATIA Series}, 2010-968, v+92 pp. 
\bibitem{klaz_III} M. Klazar, Analytic and combinatorial number theory II, {\em KAM--DIMATIA Series}, 2010-969, iv+46 pp. 
\bibitem{more_wags} C.\,J. Moreno and S.\,S. Wagstaff, Jr., {\em Sums of Squares of Integers}, 
Chapman $\&$ Hall/CRC, Boca Raton, FL, USA, 2006. 
\bibitem{niel} M.\,A. Nielsen, Introduction to the Polymath project and ``Density Hales--Jewett and Moser numbers'',
pp. 651--657 in \cite{bara_soly}.
\bibitem{poly_arch} D.\,H.\,J. Polymath,  A new proof of the density Hales--Jewett theorem, available as arXiv: 0910.3926, 
2010.
\bibitem{poly} D.\,H.\,J. Polymath, A new proof of the density Hales--Jewett theorem, {\em Ann. Math.}
{\bf 175} (2012), 1283--1327.
\bibitem{sper} E. Sperner, Ein Satz \"uber Untermengen einer endlichen Menge, {\em Math. Z.} {\bf 27} (1928),
544--548.
\bibitem{szem} E. Szemer\'edi, On sets of integers containing no $k$ elements in arithmetic progression, 
{\em Acta Arith.} {\bf 27} (1975), 199--245.
\bibitem{tao_writeup} T. Tao, Szemer\'edi's proof of Szemer\'edi's theorem, available from:\\
{\tt http://www.math.ucla.edu/\~{}tao/preprints/acnt.html}
\bibitem{tao_szem} T. Tao, A quantitative ergodic theory proof of Szemer\'edi's theorem, 
{\em Electron. J. Combin.} {\bf 13} (2006), Research Paper 99, 49 pp. 
\bibitem{tows} H. Towsner, A model theoretic proof of Szemer\'edi's theorem, available as arXiv: 1002.4456, 2011.
\bibitem{vdWa} B.\,L. van der Waerden, Beweis einer Baudetschen Vermutung, {\em Nieuw Archief voor Wiskunde}
{\bf 15} (1927), 212--216.
\end{thebibliography}
\end{document}